\documentclass[11pt,twoside]{article}

\usepackage{graphicx}
\usepackage{multirow}
\usepackage{fancyhdr}
\usepackage{array}
\usepackage{amsfonts}
\usepackage{amsmath}
\usepackage{amssymb}
\usepackage{amsthm}

\usepackage{subfigure}

\usepackage{graphics} 
\usepackage{epsfig} 
\usepackage{latexsym}

\usepackage[textwidth=16cm,textheight=22cm,top=3cm,bottom=3cm]{geometry}
\usepackage{fancyhdr}

\newtheorem{theorem}{Theorem}

\newtheorem{remark}{Remark}
\theoremstyle{definition}

\usepackage{color}

\graphicspath{{./figures/}}

\newcommand{\drop}[1]{}
\newcommand{\no}{\noindent}
\newcommand{\fer}[1]{(\ref{#1})}
\newcommand{\qtext}[1]{\quad\text{#1}}
\newcommand{\qtextq}[1]{\quad\text{#1}\quad}

\newcommand{\be}{\mathbf{e}}
\newcommand{\bF}{\mathbf{F}}

\newcommand{\bG}{\mathbf{G}}

\newcommand{\bk}{\mathbf{k}}

\newcommand{\bu}{\mathbf{u}}

\newcommand{\bw}{\mathbf{w}}
\newcommand{\bx}{\mathbf{x}}
\newcommand{\by}{\mathbf{y}}

\newcommand{\bcero}{\mathbf{0}}

\def\brho{\mbox{\boldmath$\rho$}}
\def\bpsi{\mbox{\boldmath$\psi$}}
\def\boldeta{\mbox{\boldmath$\eta$}}

\newcommand{\cL}{\mathcal{L}}
\newcommand{\cM}{\mathcal{M}}
\newcommand{\cN}{\mathcal{N}}
\newcommand{\cQ}{\mathcal{Q}}

\newcommand{\cS}{\mathcal{S}}

\newcommand{\eps}{\varepsilon}

\newcommand{\grad}{\nabla}

\newcommand{\p}{\partial}

\newcommand{\N}{\mathbb{N}}
\newcommand{\R}{\mathbb{R}}

\def\O{\Omega}

\newcommand{\abs}[1]{| #1 |}
\newcommand{\nor}[1]{\| #1 \|}

\DeclareMathOperator{\Div}{div}

\DeclareMathOperator{\tr}{tr}
\DeclareMathOperator{\real}{Re}

\title{ An economic cross-diffusion mutualistic model for cities emergence
\thanks{First author supported by the Spanish MEC Project ECO2016-76818. Second author supported by the Spanish MEC Project MTM2017-87162-P.}}

\author{Gonzalo F. de-C\'ordoba \thanks{Dpto.\ de Teor\'{\i}a e Historia Econ\'omica, Universidad de M\'alaga, 29013-M\'alaga, Spain ({\tt gfdc@uma.es})}
    \and Gonzalo Galiano \thanks{Dpt. of Mathematics, Universidad de Oviedo, c/ Calvo Sotelo, 33007-Oviedo, Spain ({\tt galiano@uniovi.es})}\footnotemark[2] }

\date{}

\begin{document}

\thispagestyle{plain}

\maketitle

\begin{abstract}
  We study an evolution cross-diffusion problem with mutualistic Lotka-Volterra reaction term to modelize the long-term spatial distribution of labor and capital. The mutualistic behavior is deduced from the gradient flow associated to profits maximization. We perform a 
linear and weakly nonlinear stability analysis and find conditions under which 
the uniform profits optimum becomes unstable, leading to pattern formation.
The patterns alternate regions of high and low concentrations of labor and capital, which may be interpreted as cities. Finally, numerical simulations based on the weakly nonlinear analysis, as well as in a finite element approximation,
are provided. 

\no\emph{Keywords: }
Cross-diffusion, mutualism, Turing instability, weakly nonlinear analysis, labor, capital, city.


\end{abstract}

\section{Introduction}
The first human settlements that can be regarded as cities had one element in common: 
the presence of a natural formation that endowed the location with a characteristic making small amounts of capital particularly productive. This is, for instance, the case of the alluvial plains of Lower Mesopotamia, where a natural levee made it very easy to irrigate large extensions of land with few tools. Moreover, the recent invention of the plow led to an extraordinary gain in productivity.

In this article, we study a mathematical model to explore the type of relation between labor and capital that leads to the formation of cities, understood as accumulation of labor and capital. Unlike other models, where this interaction is assumed to take place in a discrete spatial domain formed by patches of employment  \cite{Aly2012}, we  consider a continuous time-space  model where the relations between labor and capital leading to global growth or decay are derived solely from a production function, and where the spatial relocation is induced by attraction-repulsion mechanisms among the units of labor and capital and, secondarily, by random motion.

Current economic geography models like those introduced by Krugman \cite{Krugman1991}
or by Tabuchi and Thisse \cite{Tabuchi2010} are focused in modelling industry concentration or in explaining the emergence of \emph{central places}.   These models are set in a discrete spatial domain, being the main forces  explaining concentration  the decreasing transportation costs and the increasing returns to scale in the production of manufactured goods. A production function is said to have increasing returns to scale when the increase of one unit of all the inputs, induce an increase of output larger than one. See \cite{Layson2014} for an extended discussion about returns to scale in a production function similar to ours.  

In the model we introduce, cities emergence is a possible outcome of the action of general interacting forces that relocate labor and capital, and which is tidily connected to growth. Thus, we do not need to specify the actual forces involved in the process, that can include those of Krugman or Tabuchi and Thisse.
Furthermore, we do not make use of the assumption of increasing return to scale in the production of goods.

Our model is related to that introduced by Volpert, Petrovskii, and Zincenko \cite{Volpert2017, Zincenko2018}
for the interaction of human migration (labor, in our model) and wealth distribution (capital). Both models are formulated in terms of reaction-diffusion partial differential equations in which the reaction term establishes the growth-decline relationship between humans and capital, and where the diffusion term stands for the spatial relocation of both densities.
 However, there are important differences between both models, being perhaps the most important the role played by the cross diffusion mechanism, which is weak in Volpert's model while crucial in ours. Besides, Volpert's analysis is focused in the existence of travelling fronts, this is, in the existence of particular solutions of the problem, while our focus is in Turing's bifurcation analysis.

\subsection{A historical example} 
Before moving into the description of the mathematical model, let us motivate it by examining the economic events that happened five thousand years ago in the levees of Mesopotamia \cite{Emberling}. 

Previously to the arising of the first civilizations, the Humanity, for hundreds of thousands of years, was nearly \emph{uniformly} distributed in geographical areas where the environmental conditions were favorable. 
This uniform state suffered a perturbation caused by the Neolithic revolution. In the levees of Mesopotamia,  people concentrated along the natural ridge because it was easy to irrigate, and with the use of plows and animal force, people found their labor very productive. There was \emph{growth}. 

The increase of productivity \emph{attracted} more people who \emph{relocated} in longer stretches of the ridge producing more and more until the ridge was crowded. Because, the ridge was on limited supply, more workers and more capital faced decreasing returns to scale due to \emph{competence}. 

At some moment, the extra capital and the extra workers were \emph{repelled} by the overcrowding, finding themselves more productive somewhere else. Thus, they relocated to a, possibly, less well endowed place than the ridge.    

This process continued along the time, with the formation of new capital and labor agglomerations that led to the distribution that we observe nowadays, far from the initial prehistorical uniform equilibrium.

\subsection{The mathematical model}
After the work of Turing \cite{Turing1952},  it is well-known that the introduction of spatial mechanisms through diffusion terms in reaction systems of differential equations may change the stability of uniform equilibria, developing new ones with a non-uniform pattern-like spatial profile.

Diffusion terms allow to modelize random relocation as well as redistribution due to attractive or repulsive pressures. An important example of reaction-diffusion models is the chemotaxis model, where bacteria relocates due to attraction forces caused by a chemical stimulus. Starting with the works of Keller and Segel \cite{KellerSegel1971}, this model received much attention due to its modelling simplicity, this is, the easy identification of the biological mechanisms implied by each term of the equations; its analytical and numerical tractability, which allows to prove key mathematical properties like the existence, uniqueness and nonnegativity of the solutions; and, finally, the ability of the model to capture essential aspects of the population behavior \cite{HillenPanter2009}.

The chemotaxis model is a particular example of reaction cross-diffusion models \cite{Shigesada1979, Amann1989,  Galiano2003, Chen2004, Desvillettes2014, Galiano2014, Jungel2015, Jungel2017, Ahmed2002}.  The reaction term  induces the global growth or decline of the substances, while the diffusive term redistributes them in space. Correspondingly, our model is constructed in two steps.

\begin{enumerate}
 \item For the reaction terms, we investigated which functional form should be considered for the labor-capital relationship. As shown in \fer{eq:u0}-\fer{eq:v0}, we found that this relation is mutualistic, in contrast with previous literature, where a predator-pray model was used  \cite{BalazsiKiss2015}. 
 We delay the justification of this assumption to Section~\ref{sec:parameters}.
 
 \item For the diffusion part of the system, we assume the following general qualitative laws:
 \begin{enumerate}
 \item Labor and capital relocate randomly moving from higher to lower density regions. 
 \item Labor is repelled by high labor density regions.
 \item Labor and capital are attracted by high capital density regions. However, the attraction felt by capital is limited, due to increasing costs of relocation to such regions.
\end{enumerate}
\end{enumerate}

With these laws operating, the mathematical problem is formulated as follows. 
Let $L$ and $K$ denote the labor and capital concentrations. 
Find $L,K:[0,T]\times\bar\O\subset\R^n\to\R_+$, with $n\in\N$, such that  
\begin{align}
&  \p_t L - \Div \big( (c_1 +a_{11} L) \grad L -a_{12}L \grad K\big) = L(\alpha_1 -\beta_{11}L + \beta_{12}K) & &\text{in } Q_T, & \label{eq:u0}\\
&  \p_t K -  \Div\Big( \big(c_2 -a_{22} \frac{K}{K_s^2+K^2} \big) \grad K \Big) = K (\alpha_2 +\beta_{21}L - \beta_{22}K) & &\text{in } Q_T ,\label{eq:v0}&\\
& \grad L\cdot \nu = \grad K \cdot \nu = 0 & &\text{on }  \Gamma_T, & \label{eq:bd0}\\
& L(0,\cdot) = L_0, \quad K(0,\cdot) = K_0 & &\text{in } \O, & \label{eq:id0}
\end{align}
where $Q_T=(0,T)\times \O$, $\Gamma_T=(0,T)\times \p\O$, and the coefficients are non-negative constants. 

For the diffusion part of \fer{eq:u0}-\fer{eq:v0}, we adopted the simplest mathematical functional forms for random diffusion (terms containing $c_1,c_2$) and for the intra-repulsion and inter-attraction experienced  by labor (terms containing $a_{11}$ and $a_{12}$). However for the intra-attraction experienced by capital (term containing $a_{22}$), we incorporate a saturation function  which results on a more complex  nonlinear interaction. The reason is the following: high density capital regions are attractive for labor and capital due to work demand and capital synergies, respectively. However, this trend has a limit for the capital due to rising prices of access to real estate and labor, increment of taxes, etc. In this limit, that we capture through the saturation constant $K_s$, the attraction function reaches a global maximum at the given geographical region. Thus,  injecting additional capital in the place only makes it less attractive. See Remark~\ref{rem:saturation} for further explanations on the saturation term.

In relation to the boundary conditions \fer{eq:bd0}, we assume that labor and capital are isolated within the spatial domain $\O$, that is, there are no labor or capital flows through the boundary of $\O$. Concerning the initial data \fer{eq:id0}, we assume that $L_0$ and $K_0$ are smooth non-negative given functions.

The rest of the paper is organized as follows. In Section~\ref{sec:parameters}, we investigate the kind of reaction term which should be adopted in the model. Starting from a general production function, we show that the gradient flow of the corresponding profits function may be approximated by a Lotka-Volterra system. 
We show that standard economic assumptions imply that this system is of mutualistic type. We illustrate this  fact by considering  the  CES production function, see \fer{def:ces}.

In Section~\ref{sec:turing}, we perform a Turing stability analysis to find conditions under which the uniform equilibrium becomes unstable and a new pattern-like equilibrium arises. In geographical economy terms, we check that the spatially uniform coexistence equilibrium, which was prevalent during the long early period of human development, becomes unstable and produces a new non-uniform equilibrium in which geographical patterns arise. These geographical patterns, accumulations of labor and capital, may be interpreted as cities.

Although the linear stability theory is a useful step for understanding pattern formation, it only gives a rough indication of the patterns we should expect. Thus, in the Appendix, we perform a weakly nonlinear analysis based on the method of multiple scales that, in addition to the unstable wave-numbers provided by the linear theory, allows also to approximate the amplitude of such instabilities.

Finally, in Section~\ref{sec:numerics}, we illustrate the theoretical results by performing several numerical experiments in which the data lead to the instability of the uniform equilibrium. We approximate the full nonlinear system 
\fer{eq:u0}-\fer{eq:id0} with a finite elements scheme and, when possible, compare this approximation with that provided by the weakly nonlinear analysis.

\section{The reaction term: parameters of competence and mutualism }\label{sec:parameters}

A production function, $Y(K,L)$, is an abstraction that summarizes a technical relation between the inputs of labor and capital used in a productive process and the resulting quantity of output. 

Firms, in order to operate a production function, must pay the use of productive factors. As a compensation, they obtain income through sales, but at the cost of remunerating the worker hours used and the capital hired. The difference between the revenue sales and the factor costs, constitutes the profits:
\begin{equation*}
  \Pi(K,L)=pY(K,L)-wL-rK, 
\end{equation*}
where $p, ~w,~r$ are functions of time describing, respectively, the sale price of one unit of output, the wages paid for a unit of work, and the rental rate of capital used in production. Function $p$ can be taken as \textit{numeraire}, normalizing it to $1$.

The optimal trajectory $(K(t), L(t))$ is determined by the maximization of profits. 
The first optimality condition is 
\begin{align}
\label{eq:equiY}
  \p_L Y(K ,L ) -w = 0, \quad   \p_K Y(K ,L ) -r =  0. 
\end{align}
To compute the optimal trajectory, one must prescribe wages and rental rates or, otherwise, state suitable relationships among $L,~K,~w$ and $r$ that close the system of equations. 
A simplifying assumption consists on supposing that $w$ and $r$ are close to some stable constant values, say $w^*, ~r^*$. Then,  the profits optimum for these values, $(L_e^*, K_e^*)$, is obtained from the optimality conditions \fer{eq:equiY}.  

A dynamics approximating this optimal value is given by the nonlinear problem corresponding to the gradient flow $\p_t(K, L)^t = \grad \Pi(K, L)$. In a further simplification, we may linearize the gradient flow by performing the substitution 
$ \grad \Pi(K, L) \approx \frac{1}{2}H(\Pi_e^*) (L-L_e^*,  K-K_e^*)^t$,  
where $H(\Pi_e^*)$ is the Hessian matrix of $\Pi$ at $(L_e^*, K_e^*)$.

In order to accommodate the gradient flow to the Lotka-Volterra terms, we use a modified gradient flow, see Remark~\ref{rem:flow}, that leads to the system 
\begin{align}
& \p_t L = L \Big( \tilde\alpha_1 + \tilde\beta_{11}L + \tilde\beta_{12}K \Big) ,\label{eq:LVeq1}\\
& \p_t K = K \Big( \tilde\alpha_2 + \tilde\beta_{21}L + \tilde\beta_{22}K \Big) ,
\label{eq:LVeq2}
\end{align}
with 
\begin{align}
\label{def:alphabeta}
\tilde B:=(\tilde\beta_{ij})=H(\Pi_e^*),\quad (\tilde\alpha_1,\tilde\alpha_2)^t = -\tilde B ( L_e^* , K_e^*)^t. 
\end{align}
 In particular, notice that the profits optimum $(L_e^* , K_e^*)$ coincides with the non-trivial equilibrium of the Lotka-Volterra system. 
In addition, $\tilde\beta_{12}=\tilde\beta_{21}$ and the second optimality condition for the profits maximization, i.e. that $H(\Pi_e^*)$ is positive definite, imposes the necessary condition $\det(\tilde B) >0$.

Since, by definition, $H(Y_e^*) = H(\Pi_e^*)$, the signs of $\tilde\beta_{ij}$ are determined by the production function. In economic terms, 
the Law of \textit{decreasing marginal productivities} states that the more one factor is used, the smaller is the increase in output obtained. This implies that $\p_{LL}Y$ and $\p_{KK}Y$ are negative. However, the crossed marginal productivities are increasing, i.e., $\p_{LK}Y=\p_{KL}Y >0$. The reason is that if, for instance, labor is abundant, an increase in capital will increase the productivity of workers, and as a consequence it will improve the output, as the farmer with naked hands is much less productive than the farmer with a plow. 
Therefore, one expects the following signs for the quadratic coefficients of the Lotka-Volterra term:
\begin{align*}
\tilde\beta_{11} <0,\quad \tilde\beta_{22}<0,\qtextq{and }\tilde\beta_{12}=\tilde\beta_{21}>0,                                                                                                            \end{align*}
capturing a mutualistic model. 
\begin{remark}
\label{rem:flow}
Suppose that the dynamical system $\p_t \bu = \bF(\bu)$, for some differentiable function  $\bF:\R^2\to\R^2$, has a nontrivial asymptotically stable equilibrium $\bu^*$ with non-negative components. Let $J\bF(\bu^*)$ denote its Jacobian matrix at $\bu^*$, so that 
$\tr(J\bF(\bu^*))<0$ and $\det(J\bF(\bu^*)) >0$. 

Consider the system $\p_t \bw =  \bG(\bw)$, with $\bG_i(\bw) = w_i \bF_i(\bw)$, for $i=1,2$, having the equilibria $\bcero$ and $\bu^*$.  We have 
\begin{align*}
\tr(J\bG(\bw)) = \sum_{i= 1}^2w_i \p_i\bF_i(\bw), \quad
\det(J\bG(\bw)) = w_1w_2 \det(J\bF(\bw)),
\end{align*}
and, consequently, $\bu^*$ is its only stable equilibrium.
\end{remark}

\subsection{Example: The CES production function}
The production function of \emph{constant elasticity of substitution} (CES) \cite{Layson2014} has the following functional form:
\begin{equation}
  Y(K,L) = A\left (\alpha K^{\eta} + \beta L^{\eta}\right)^{\frac{\epsilon}{\eta}}, \label{def:ces}
\end{equation}
where $A$ is the total factor productivity, a non-negative function of time 
which is a scaling parameter used to accommodate the measures of $K$ and $L$ in relation to $Y$, and where $\alpha, \beta$ are the income share parameters. We assume the usual parameters range
\begin{align}
 \label{param}
\alpha,\beta,\epsilon \in (0,1),\quad  \eta < 1, \quad \eta\neq0, \qtext{and } \epsilon \geq \eta.
\end{align}
For this choice of the production function, we have 
\begin{align}
&  \p_L Y(K ,L ) = A \epsilon[ \alpha K ^{\eta}+\beta L ^{\eta}]^{\frac{\epsilon}{\eta}-1} \beta L ^{\eta-1}  , \label{eq:ceslaborfoc2}\\
&  \p_K Y(K ,L ) = A\epsilon[ \alpha K ^{\eta}+\beta L ^{\eta}]^{\frac{\epsilon}{\eta}-1} \alpha K ^{\eta-1} , \label{eq:cescapitalfoc2}
\end{align}
and 
\begin{align}
\label{def:hessian}
 H(\Pi) = \begin{pmatrix}
           \dfrac{\beta(\epsilon-\eta)Q^{-1}wL^{\eta} -(1-\eta)w }{L} & \alpha(\epsilon-\eta)Q^{-1}w K^{\eta-1}\\
           \beta(\epsilon-\eta)Q^{-1}r L^{\eta-1} & \dfrac{\alpha(\epsilon-\eta)Q^{-1}rK^{\eta} - (1-\eta)r}{K}
          \end{pmatrix},
\end{align}
where $Q = \alpha K^{\eta} +\beta L^{\eta}$. 

In the following theorem, we show that for a CES production function the signs of the reaction part of the system \fer{eq:LVeq1}-\fer{eq:LVeq2} 
 are those of  a mutualistic Lotka-Volterra model. This gives support to our assumption on the form of the reaction part of the system \fer{eq:u0}-\fer{eq:v0}.
We also check that the equilibrium of this model coincides with the optimum of the profits function.
\begin{theorem}
\label{proposition:1}
 Let $w^*,~r^*$ be positive numbers, and consider the CES production function \fer{def:ces} with parameters satisfying \fer{param}. Then, the optimum of the profits function, $\Pi$,  is $( L_e^*,K_e^{*})$, with 
\begin{align}
& (L_e^*)^{1-\eps}=  \frac{ \rho }{w^*}\beta^{\epsilon/\eta}, \quad  
(K_e^*)^{1-\eps}=  \frac{ \rho }{r^*}\alpha^{\epsilon/\eta}, \quad  
  \rho  = 
A\epsilon \Big( 1+ \big(\frac{w^* }{r^* }\big)^{\frac{\eta}{1-\eta}}\big(\frac{\alpha}{\beta}\big)^{\frac{1}{1-\eta}} \Big)^{\frac{\epsilon-\eta}{\eta}}. \nonumber
\end{align}
Moreover, $(K_e^*, L_e^*)$ is also the coexistence equilibrium of the Lotka-Volterra system \fer{eq:LVeq1}-\fer{eq:LVeq2} with coefficients given by \fer{def:alphabeta}, which satisfy
\begin{align*}
\tilde\alpha_1 >0,\quad \tilde\alpha_2 >0,\quad \tilde\beta_{11}<0,\quad \tilde\beta_{22}<0,\quad \tilde\beta_{12}=\tilde\beta_{21}\geq0,\quad \det(\tilde B) >0.
\end{align*}
In addition, if $\epsilon >\eta$ then $\tilde\beta_{12}=\tilde\beta_{21}>0$.
\end{theorem}
\no\emph{Proof. }
Dividing the first optimality condition by the second, see \fer{eq:equiY}, and using the explicit expressions \fer{eq:ceslaborfoc2}-\fer{eq:cescapitalfoc2} yields
\begin{equation*}
K_e^*=\left(\frac{r^* \beta}{w^* \alpha}\right)^{\frac{1}{\eta-1}}L_e^*.
\end{equation*}
Replacing in \fer{eq:equiY} and rearranging terms  we obtain the expression of $( L_e^*,K_e^{*})$ stated in the proposition. 
The optimum $(L_e^*,K_e^*)$ is also the coexistence equilibrium of \fer{eq:LVeq1}-\fer{eq:LVeq2} by construction, see \fer{def:alphabeta}. 
From the definition of the Hessian, see \fer{def:hessian}, and the condition $\epsilon\geq\eta$ (resp. $\epsilon >\eta$) we easily deduce that $\beta_{12}=\beta_{21}\geq0$ (resp. $\beta_{12}=\beta_{21}>0$). Moreover, after some computations, we obtain 
\begin{align*}
\det(\tilde B)= (1-\epsilon)(1-\eta)\frac{w^*r^*}{L_e^*K_e^*} >0.
\end{align*}
Therefore, necessarily $\tilde\beta_{11}\tilde\beta_{22}>0$. Suppose that both numbers are positive. Then, using their definition, we deduce
\begin{align*}
 \alpha\beta\frac{(L_e^*)^\eta (K_e^*)^\eta}{(\alpha (L_e^*)^\eta+\beta (K_e^*)^\eta)^2} > \Big(\frac{1-\eta}{\epsilon -\eta}\Big)^2.
\end{align*}
Set $x = \alpha (L_e^*)^\eta$, $y =  \beta (K_e^*)^\eta$ and $\delta = (1-\eta)^2/(\epsilon-\eta)^2$. The above inequality translates to 
$ \delta (x^2+y^2) +(2\delta-1)x y <0$, 
which is not possible if $2\delta >1$. And this is indeed the case, due to the conditions $\eps<1$ and $\eta< 1$. The contradiction arises after the assumption $\tilde\beta_{11}>0$ and $\tilde\beta_{22}>0$. Therefore, these quantities must be  both negative, since their product is positive. Finally, from the definition of $\tilde\alpha_i$, see \fer{def:alphabeta},
we get $\tilde\alpha_1 = w^* (1-\epsilon)>0$, and $\tilde\alpha_2 = r^* (1-\epsilon)>0$.
$\Box$

\section{Turing instability analysis}\label{sec:turing}

Before starting the stability analysis, we reduce slightly the profusion of parameters in our model by introducing the assumption that the mutualistic coefficients $\beta_{12}$, $\beta_{21}$ of the reaction terms are positive. This is, for instance, the case for a CES production function with $\epsilon >\eta$, see Theorem~\ref{proposition:1}. Thus, we perform the  following change of variables and parameters:
\begin{align*}
& \bar L = \beta_{21} L,\quad \bar K = \beta_{12}K, \quad \bar L_0 = \beta_{21}L_0, \quad \bar K_0 = \beta_{12}K_0, \quad  \bar K_s = \beta_{12}K_s,\\
& a_{1} = \frac{a_{11}}{\beta_{21}}, \quad a_{2} = a_{22}\beta_{12}, \quad b = \frac{a_{12}}{\beta_{12}},
 \quad \beta_{1} = \frac{\beta_{11}}{\beta_{21}},  \quad \beta_{2} = \frac{\beta_{22}}{\beta_{12}},\\
& \bar t = \frac{t}{\gamma},\quad \bar x = \frac{x}{\sqrt{\gamma}},
\end{align*}
for some time-space re-scaling factor $\gamma >0$.  
This change renders problem \fer{eq:u0}-\fer{eq:id0} to the following form (omitting bars)
\begin{align}
&  \p_t L - \p_x \big( (c_1 +a_{1} L) \p_x L -b L \p_x K\big) = \gamma L(\alpha_1 -\beta_{1}L + K) & &\text{in } Q_T, & \label{eq:u}\\
&  \p_t K -  \p_x\Big( \big(c_2 -a_{2} g(K) \big) \p_x K \Big) = \gamma K (\alpha_2 +L - \beta_{2}K) & &\text{in } Q_T ,\label{eq:v}&\\
& \p_x L\cdot \nu = \p_x K \cdot \nu = 0 & &\text{on }  \Gamma_T, & \label{eq:bd}\\
& L(0,\cdot) = L_0, \quad K(0,\cdot) = K_0 & &\text{in } \O, & \label{eq:id}
\end{align}
where $(0,T)$ and $\O$ are redefined according to the scaling factor $\gamma$. We collect here some assumptions and necessary conditions for the well-posedness of the problem:
\begin{enumerate}
 \item The saturation function $g \in C^2(\R_+)$ is non-negative and reaches its maximum at $K_s>0$. 
 
 \item The coefficients $c_i$, $b$, and $a_i$, for $i=1,2$, are non-negative. In particular, we assume the ellipticity conditions
 \begin{align*}
  c_1+a_1 >0,\quad c_2 > a_2 g(K_s).
 \end{align*}
 \item In agreement with Theorem~\ref{proposition:1}, the coefficients $\alpha_i$ and $\beta_i$, for $i=1,2$, are non-negative. Moreover, dominance of competitiveness over mutualism, captured by $\beta_1\beta_2 >1$, is also assumed.
\end{enumerate}
Notice that after the change of unknowns  leading to the system  \fer{eq:u}-\fer{eq:id}, the coexistence equilibrium is expressed as
\begin{align}
\label{def:eq}
(L^*, K^*) = \Big(\frac{\alpha_2+\alpha_1\beta_2}{\beta_1\beta_2-1}, \frac{\alpha_1+  \alpha_2\beta_1}{\beta_1\beta_2-1}\Big).
\end{align}

In the following theorem we establish that diffusion induced instability occurs when 
the attraction felt by labor toward regions with a high density of capital is large enough.
\begin{theorem}
The coexistence equilibrium \fer{def:eq}, which is stable for the reaction system associated to \fer{eq:u}-\fer{eq:id}, becomes unstable for the whole reaction-diffusion system \fer{eq:u}-\fer{eq:id} when $b>b_c$, see \fer{def:bc}. This instability leads to pattern formation if $\gamma$ is large enough.  
\end{theorem}
\no\emph{Proof. }
Let $(\tilde L, \tilde K)$ be a small perturbation of this equilibrium satisfying \fer{eq:u}-\fer{eq:bd} and set 
$(w_1,w_2) = (\tilde L -L^*, \tilde K-K^*)$. The linearized problem for $\bw=(w_1,w_2)$ is given in matrix form by
\begin{align}
&  \p_t \bw -  Q \Delta \bw  = \gamma R\bw & &\text{in } Q_T, & \label{eq:ul}\\
& \grad u_1\cdot \nu = \grad u_2 \cdot \nu = 0 & &\text{on }  \Gamma_T, & \label{eq:bdl}\\
& u_1(0,\cdot) = \tilde L(0,\cdot)-L_0, \quad u_2(0,\cdot) = \tilde K(0,\cdot)-K_0 & &\text{in } \O, & \label{eq:idl}
\end{align} 
where 
\begin{align}
\label{def:Q}
R= \begin{pmatrix}
-\beta_{1}L^* & L^*\\
K^* & -\beta_{2} K^*
\end{pmatrix},
\qquad 
Q= \begin{pmatrix}
c_1 +a_{1} L^* & -bL^*\\
0 & c_2 -a_{2} g(K^*)
\end{pmatrix}.
\end{align}
Taking into account the boundary conditions, we search for particular solutions of \fer{eq:ul}-\fer{eq:idl} of the form $\bw(t,x) \propto e^{\lambda t+i\bk\cdot x}$, where $\lambda$ represents the linear growth rate and $\bk$ is the wave number of the perturbation.   Introducing this expression of $\bw$ in \fer{eq:ul}-\fer{eq:idl} we are led to the eigenvalue problem 
\begin{align}
\label{def:Ak}
A_k \bw = \lambda \bw, \qtext{for }A_k = \gamma R-k^2Q \text{ and }k=\abs{\bk}.
\end{align}
The corresponding characteristic equation 
yields the eigenvalues 
\begin{align*}
 \lambda_{jk} = \frac{1}{2}\Big( \tr(A_k) + (-1)^j \sqrt{\tr(A_k)^2 -4\det(A_k)}\Big), \quad j=1,2,
\end{align*}
so that, for Turing instability to occur,  one or more of the eigenvalues $\lambda_{jk}$ must have a positive real part.

For $k=0$ (no diffusion), we have $\tr(A_0) =  \gamma \tr (R) <0$ and $\det(A_0) =\gamma\det(R)  >0$, implying $\real(\lambda_{j0})<0$ for $j=1,2$. This is, the mode corresponding to the wave number $k=0$ is always stable. 

For $k>0$ we have $\tr(A_k) = \gamma\tr (R) - k^2 \tr(Q) <0$,
and therefore, the only way which may lead to a positive $\real(\lambda_{jk})$ involves $\det(A_k)$ being negative. We introduce  the convex quadratic polynomial 
\begin{align*}
h(k^2) := \det(A_k)=\det(Q) k^4 + \gamma  k^2 q+\gamma^2 \det(R), 
\end{align*}
with $q=m_2-b m_1$, and $m_1,m_2$ given by the positive constants 
\begin{align*}
 m_1 = L^*K^*  ,\quad 
 m_2 = \beta_{1}L^* (c_2- a_{2}g(K^*)) +  \beta_{2}K^* (c_1+a_{1}L^* ) .
\end{align*}
The minimum of $h$ is attained at 
\begin{align}
\label{def:km}
 k^2_m = -\frac{\gamma q}{2\det(Q)},
\end{align}
which is a real root if $q<0$. Hence, a necessary condition for instability is 
\begin{align}
\label{cond:nec}
b>m_2/m_1. 
\end{align}
The corresponding minimum value is 
\begin{align*}
 \min(h(k_m^2)) = \gamma^2\Big(\det(R) - \frac{q^2}{4 \det(Q)} \Big), 
\end{align*}
which is negative for $q < -2\big(\det(R)\det(Q)\big)^{1/2}$,
or, in  terms of the \emph{bifurcation parameter} $b$, if $b>b_c$, with the critical bifurcation value given by  
\begin{align}
\label{def:bc}
 b_c =\frac{m_2+2\big(\det(R)\det(Q)\big)^{1/2}}{m_1}.
\end{align}
Observe that, in particular, $b>b_c$ implies the necessary condition \fer{cond:nec}.
The critical wavenumber corresponding to the critical bifurcation values is found by replacing $q^c = m_2-b_cm_1$ in \fer{def:km}, yielding
\begin{align*}
 k_c^2 = \gamma \Big(\frac{\det(R)}{\det(Q)}\Big)^{1/2} = \gamma \Big(\frac{(\beta_1\beta_2-1)L^*K^*}{(c_1+a_1L^*)(c_2-a_2g(K^*))}\Big)^{1/2}.
\end{align*}
Thus, for $b>b_c$ the system has a range $(k_1^2, k_2^2)$ of unstable wave-numbers, where $k_1^2,k_2^2$ are the roots of $h$. Pattern formation will arise if the spatial spread of $\O$ is large enough so that at least one of the modes admitted by the boundary conditions lies in the interval $(k_1^2, k_2^2)$. Since the roots of $h$ depend linearly on $\gamma$, the time-space scale parameter, pattern formation  in the original model will be observed if this parameter is large enough.  
$\Box$

To gain some insight into the bifurcation condition, we replace the expressions of $m_1,~m_2,~\det(R)$ and $\det(Q)$ in \fer{def:bc} to obtain that $b>b_c$ is equivalent to 
\begin{align*}
L^*K^* b > \beta_1 L^* X^2 +\beta_2 K^* Y^2 + 2 (L^* K^*\det(B))^{1/2}XY,
\end{align*}
where $X^2=c_2-a_2g(K^*)$ and $Y^2=c_1+a_1L^*$. Completing squares, 
and using the elementary inequality $(x+y)^2\leq 2(x^2+y^2)$, we deduce the sufficient condition 
\begin{align*}
L^*K^* b + 2XY \sqrt{L^* K^*} \big(\sqrt{\beta_1\beta_2}-\sqrt{\det(B)}\big) > 2(\beta_1 L^* X^2 +\beta_2 K^* Y^2 ).
\end{align*}
Since $\det(B)=\beta_1\beta_2-1$, we obtain the following, more strict, sufficient condition, providing  a clear interpretation of the bifurcation condition in terms of the diffusion coefficients,
\begin{align}
\label{cond:bif}
 \frac{b}{2} + \frac{\beta_1}{K^*} a_2g(K^*) > \frac{\beta_1}{K^*} c_2 +  \frac{\beta_2}{L^*}c_1 + \beta_2a_1.
\end{align}
As expected, we see that the diffusion coefficients promoting stability are those representing random relocation and labor intra-repulsion, while those promoting instability are the intra-attraction coefficients of both labor and capital. 

In which respect to the relationship between instability and the kinematic terms, let us assume the simplification $\alpha_1=\alpha_2 =:\alpha$. Replacing these values in the expressions of the equilibrium, and then in \fer{cond:bif}, we get the sufficient condition 
\begin{align}
\label{cond:bif2}
 \frac{b}{2} + \frac{\det(B)}{\alpha} a_2g(K^*) > \frac{\det(B)}{\alpha} \big(c_1 +  c_2)+\beta_2 a_1,
\end{align}
from where we interpret that if there is a high increase of labor and capital, expressed by a small intra-competitive behavior and a large intrinsic growth,  instability will arise.

\section{Numerical simulations}\label{sec:numerics}
We approximate the solutions of the nonlinear problem \fer{eq:u}-\fer{eq:id} by  employing two methods: the finite element method (FEM) and the weakly nonlinear approximation (WNL). The latter is explained in detail in the Appendix, but we give here a short summary of the main ideas. See also \cite{Gambino2012, Gambino2013}.

Let $\eps$ be a small control parameter representing the dimensionless distance from the critical threshold $\eps^2 = (b-b_c)/b_c$. The idea of the weakly nonlinear analysis 
is to expand the perturbation, $\bw$, the bifurcation parameter, $b$, and the time scale, $t$, in powers of $\eps$, 
\begin{align}
 &b = b_c + \eps b_1 +\eps^2 b_2 +\eps^3 b_ 3 + \cdots \label{expansion1},\\
 &\bw = \eps \bw_1 +\eps^2 \bw_2 +\eps^3 \bw_ 3 + \cdots \label{expansion2}\\
 & \p_t = \eps \p_{T_1} +\eps^2 \p_{T_2} + \eps^3 \p_{T_3} + \cdots \label{expansion3}
\end{align}
After substitution in the nonlinear system satisfied by $\bw$, see \fer{eq:wexp}, and equating with respect to the order of $\eps$, we obtain a chain of linear problems for $\bw_i$, for $i=1,2,\ldots$, whose solutions are of the form $\bw_i(T_1,T_2,\ldots, x) = A(T_1,T_2,\ldots) \bu_i(x)$, with $A$ denoting amplitude. In particular, the third order approximation leads to the cubic Stuart-Landau equation for the amplitude,
\begin{align*}
 \p_{T_2} A = \sigma A - \ell A^3,
\end{align*}
where the growth rate coefficient $\sigma$ is positive. The dynamics of this  equation is divided into two different cases according to the sign of the Landau constant $\ell$: the supercritical case ($\ell>0$), and the subcritical case ($\ell<0$).

In the supercritical case, the equilibrium solution $A_\infty =\sqrt{\sigma/\ell}$ is stable, and  represents the asymptotic value of the amplitude $A$. The corresponding solution, corrected to satisfy the Neumann boundary conditions, is given by
\begin{align}
\label{wnl2}
\bw =\eps\brho \sqrt{\frac{\sigma}{\ell}} \cos(\bar{k_c}x) 
  + \eps^2\frac{\sigma}{\ell}\big(\bw_{20}+
  \bw_{22}\cos(2\bar{k_c}x)\big) + O(\eps^3),
\end{align}
where $\bar{k_c}$ is the first integer or semi-integer to become unstable when $b$ passes the critical value $b_c$. 

In the subcritical case, the cubic Stuart-Landau equation does not give any valid information for the amplitude, and a higher degree expansion on powers of $\eps$ must be considered, leading to higher order Stuart-Landau equations.

For the FEM approximation, we used the open source software \texttt{deal.II} \cite{dealII83} to implement a time semi-implicit 
scheme with a spatial linear-wise finite element discretization. For the time discretization, we take 
in the experiments a uniform time partition of time step $\tau$. For the spatial  discretization, we take a uniform partition of the interval $\O =[0,2\pi]$.

Let, initially, $t=t_0=0$ and  set $(L^0, K^0)=(L_0, K_0)$. For $n\geq 1$, the discrete problem is: Find $L^{n}, K^n \in S^h$ such that  
\begin{align}
 \frac{1}{\tau}\big( L^n-L^{n-1} , \chi )^h
& + \big( (c_1+a_1 L^{n}) \grad L^n - b L^n\grad K^n  ,\grad\chi \big)^h \label{eq:disc1}\\
& =\big( \gamma L^n (\alpha_1 -\beta_1 L^n + K^n) , \chi )^h, \nonumber\\
\frac{1}{\tau}\big( K^n-K^{n-1} , \chi )^h
& + \big( (c_2-a_2g(K^{n})) \grad K^n  ,\grad\chi \big)^h \label{eq:disc2} \\
& =\big( \gamma K^n (\alpha_2 + L^n -\beta_2 K^n) , \chi )^h, \nonumber
\end{align}
for every $ \chi\in S^h $, the finite element space of piecewise $\mathbb{Q}_1$-elements. 
Here, $(\cdot,\cdot)^h$ 
stands for a discrete semi-inner product on $\mathcal{C}(\overline{\Omega} )$.

Since \fer{eq:disc1}-\fer{eq:disc2} is a nonlinear algebraic problem, we use a fixed point argument to 
approximate its solution,  $(L^n,K^n)$, at each time slice $t=t_n$, from the previous
approximation $(L^{n-1},K^{n-1})$.  Let $L^{n,0}=L^{n-1}$ and $K^{n,0}=K^{n-1}$. 
Then, for $k\geq 1$ the linear problem to solve is: Find $(L^{n,k},K^{n,k})$ such that for for all $\chi \in S^h$ 
\begin{align*}
 \frac{1}{\tau}\big( L^{n,k}-L^{n-1} , \chi )^h
& + \big( (c_1+a_1 L^{n,k-1}) \grad L^{n,k} - b L^{n,k-1}\grad K^{n,k}  ,\grad\chi \big)^h \\
& =\big( \gamma L^{n,k} (\alpha_1 -\beta_1 L^{n,k-1} + K^{n,k-1}) , \chi )^h, \\
\frac{1}{\tau}\big( K^{n,k}-K^{n-1} , \chi )^h
& + \big( (c_2-a_2g(K^{n,k-1})) \grad K^{n,k}  ,\grad\chi \big)^h  \\
& =\big( \gamma K^{n,k} (\alpha_2 + L^{n,k-1} -\beta_2 K^{n,k-1}) , \chi )^h.
\end{align*}
We use the stopping criteria 
 \begin{equation*}
 \max \big(\nor{L^{n,k}-L^{n,k-1}}_2, \nor{K^{n,k}-K^{n,k-1}}_2 \big) <\text{tol}_{FP},
 \end{equation*}
for values of $\text{tol}_{FP}$ chosen empirically, and set $(L^n, K^n)=(L^{n,k},K^{n,k})$. 
Finally, we integrate in time until a numerical stationary solution, 
$(L^S, K^S)$, is achieved. This is determined 
by 
\begin{equation*}
\max \big(\nor{L^{n,1}-L^{n-1}}_2, \nor{K^{n,1}-K^{n-1}}_2 \big) <\text{tol}_S,
\end{equation*}
where $\text{tol}_S$ is chosen empirically too. 

\begin{table}[!t]
\centering
{\small
\begin{tabular}{||l||l|l|l|l||}
\hline\hline 
 & \multicolumn{4}{|c||}{Experiments} 
 \\
 \hline  
                      & 1            & 2            & 3            & 4  \\\hline
 $w^*,~r^*$           & 1,~0.3       & 0.4,~0.3      & 0.95,~0.95  &    0.4,~0.5\\
 $\alpha,~\beta$      & 0.3,~0.6     & 0.3,~0.6     & 0.29,~0.3    &   0.3,~0.6   \\
 $\epsilon,~\eta$     & 0.5,~0.2     & 0.5,~0.2     &  0.75,~0.1   &   0.5,~0.2   \\
 $\alpha_1,~\alpha_2$ & 0.5,~0.15    & 0.2,~0.15    &  0.24,~0.24  & 0.05,~0.15   \\
 $\beta_1,~\beta_2$   & 2.35,~2.47   & 0.48,~8.38   &  0.66,~1.97  & 4.5e-2,~5.4   \\
 $c_1,~c_2$           & 0.01,~0.01   & 0.01,~0.01   &  0.1,~0.1    & 0.01,~0.01   \\
 $a_1,~a_2$           & 0.3,~3e-4    & 7e-3,~0.01   &  0.41,~2.4   & 3.72,~2.7e-5 \\
 $(L^*,K^*)$          & 0.29,~0.18   & 0.6,~0.09    &  2.35,~1.31  & 1.96,~0.03   \\
 $b_c,~k_c$           & 1.56,~3.99   & 0.42,~6      &  1.14,~4.02  & 205.9,~3.51   \\
 $T_s,~\text{MSE}$    & 1557,~1.7e-2 & 1630,~1.4e-2 & 240,~4.7e-2   & 5000   \\
 \hline\hline  
\end{tabular}

\caption{Parameters of the CES production function, resulting (scaled) Lotka-Volterra coefficients, rescaled diffusion coefficients, equilibrium, critical bifurcation and wave numbers, time to stationary state and mean square error between the FEM solution and the WNL solution.}
\label{tab:data}
}
\end{table}

\begin{figure}[t!]
\centering
{\includegraphics[width=6cm,height=4cm]{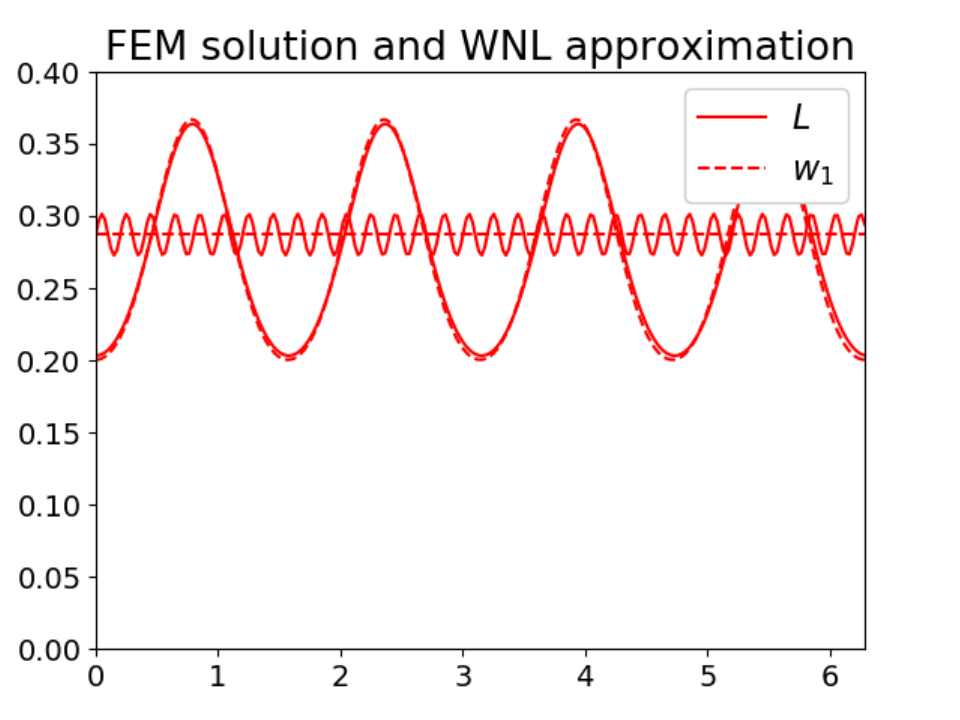}}
 \hspace{0.1cm}
{\includegraphics[width=6cm,height=4cm]{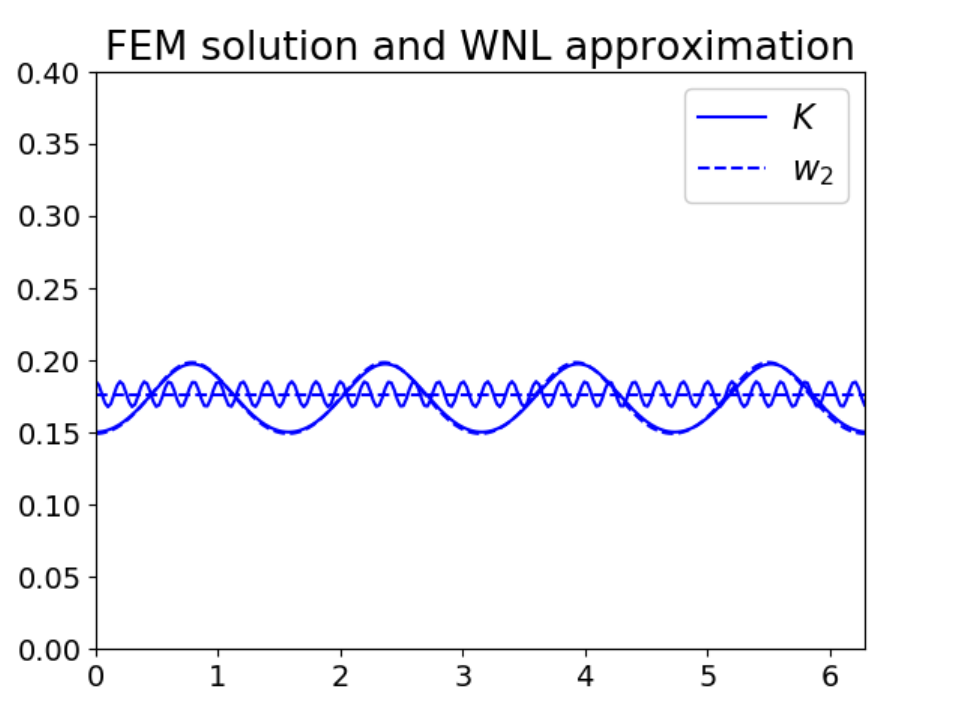}}\\
{\includegraphics[width=6cm,height=4cm]{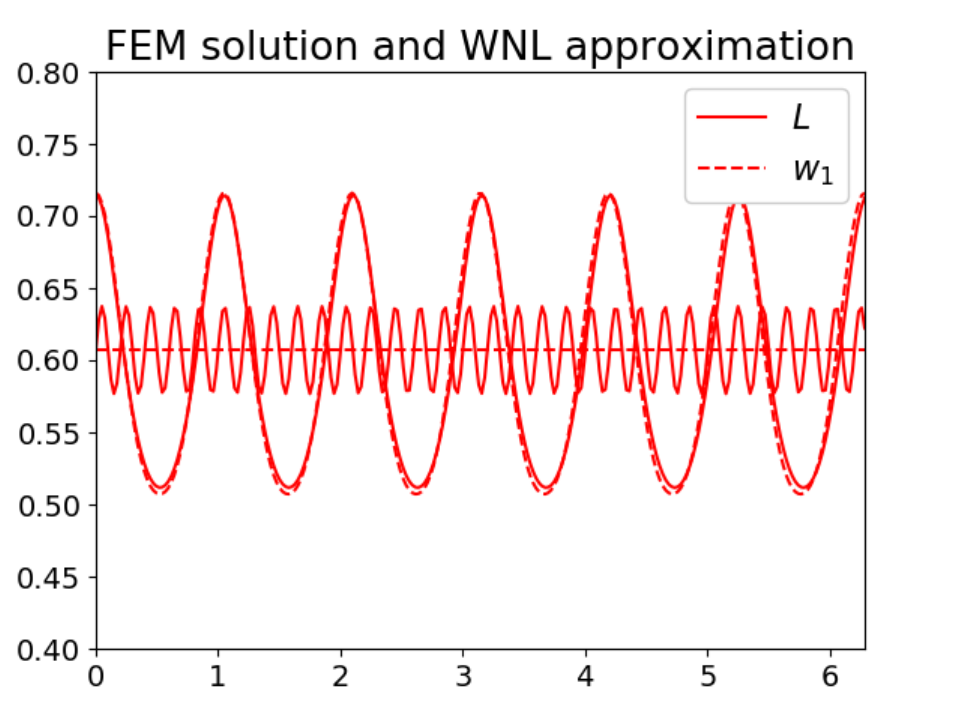}}
 \hspace{0.1cm}
{\includegraphics[width=6cm,height=4cm]{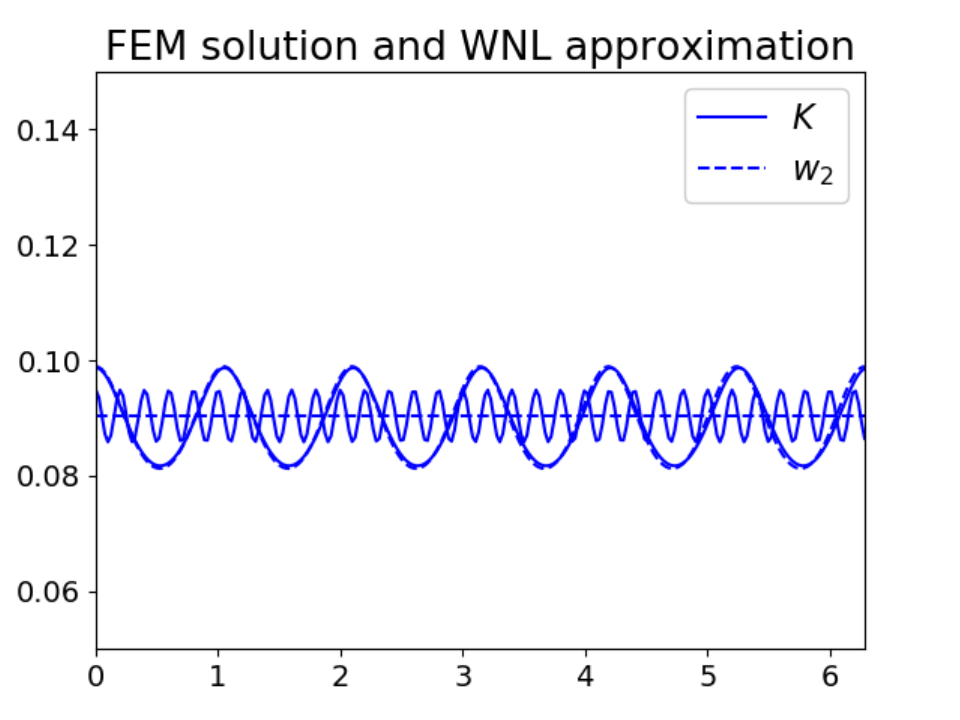}}\\
{\includegraphics[width=6cm,height=4cm]{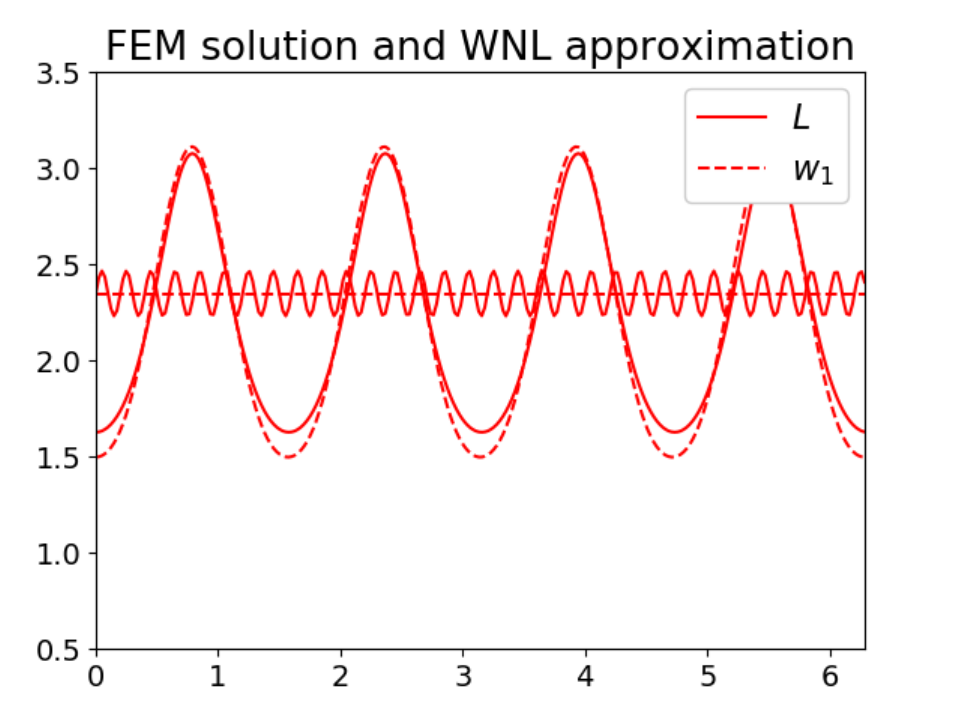}}
 \hspace{0.1cm}
{\includegraphics[width=6cm,height=4cm]{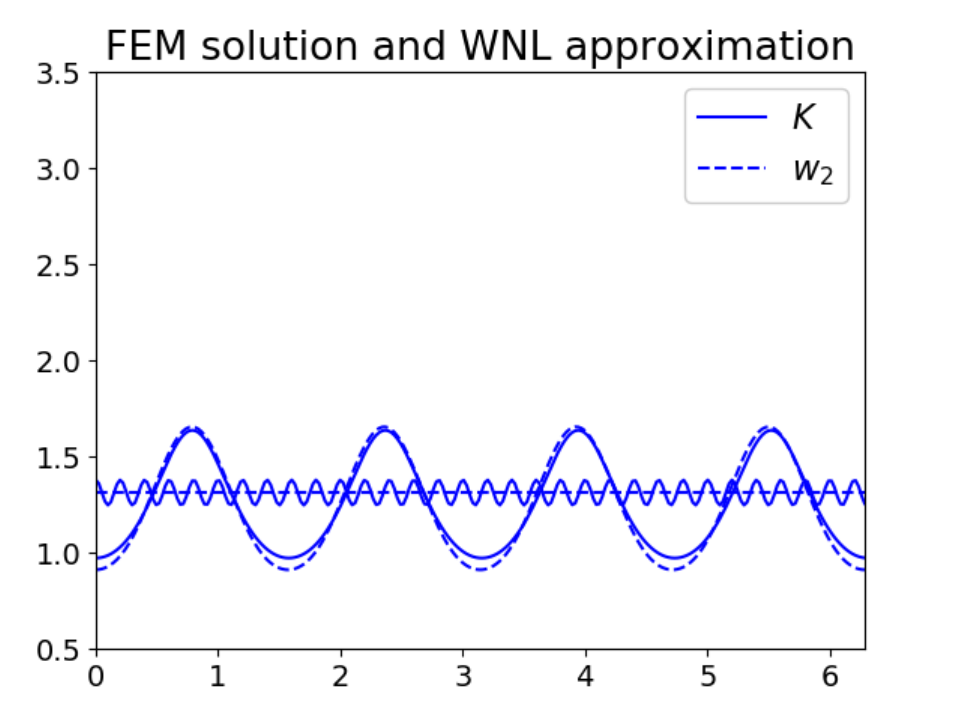}}
\caption{{\small Supercritical case. First row: Experiment 1. New equilibrium reached after the onset of instabilities due to a small perturbation around the uniform equilibrium $(L^*,K^*)\approx(0.287,0.176)$. Second row: Experiment 2.  Like previous, with  $(L^*,K^*)=(0.6,0.09)$. Mind the different ordinate scales. Third row: Experiment 3.  Like previous, with  $(L^*,K^*)\approx(2.35,1.31)$.
}} 
\label{exp1_fig}
\end{figure}

\begin{figure}[t!]
\centering
{\includegraphics[width=6cm,height=4cm]{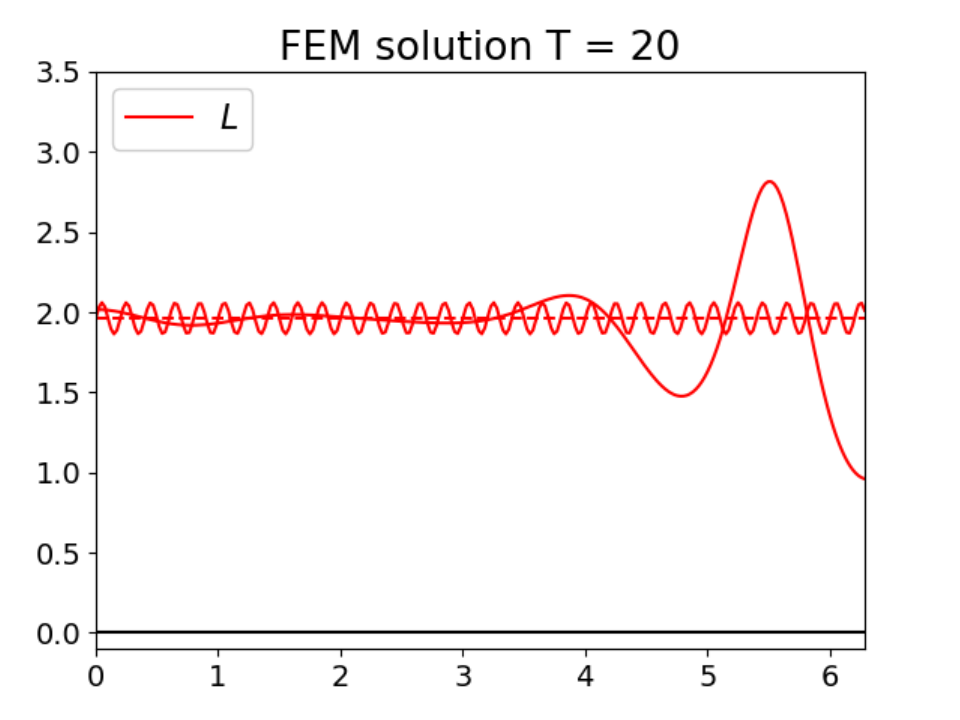}}
 \hspace{0.1cm}
{\includegraphics[width=6cm,height=4cm]{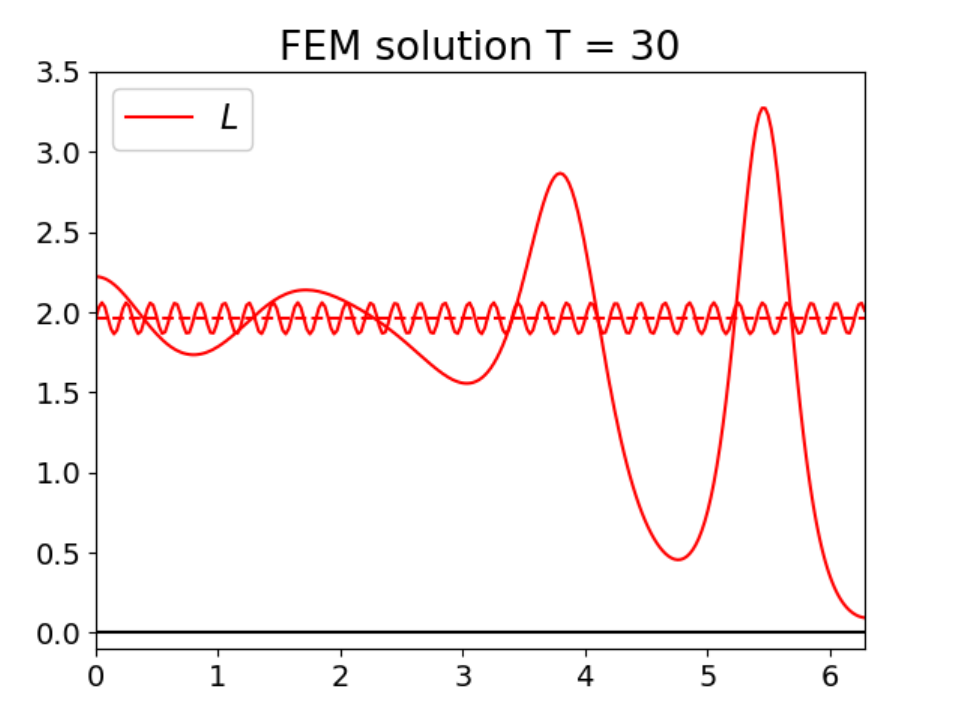}}\\
{\includegraphics[width=6cm,height=4cm]{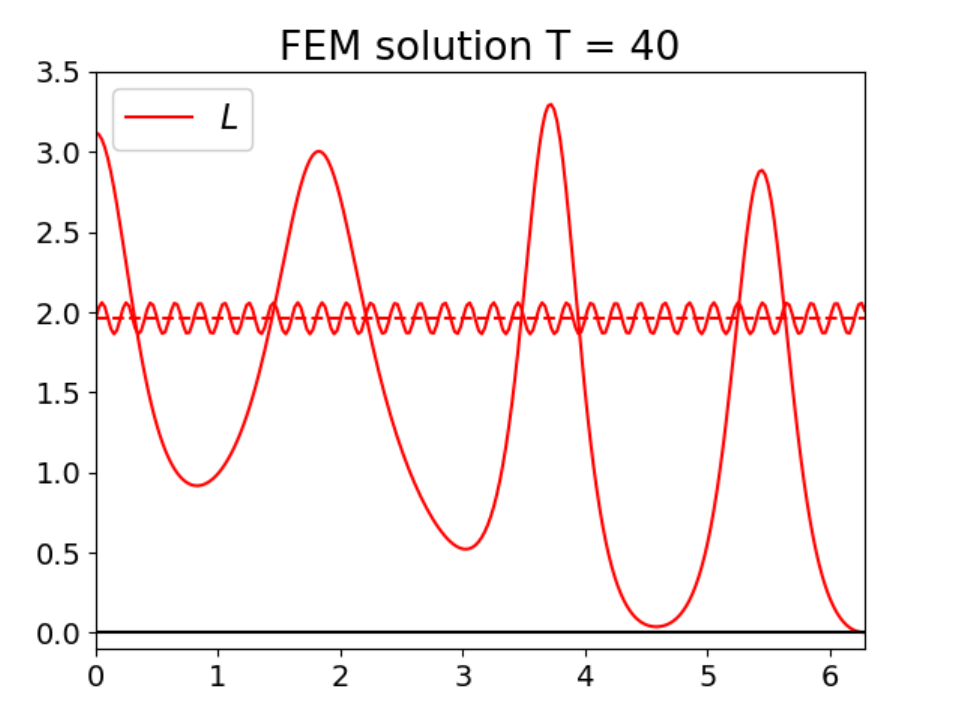}}
 \hspace{0.1cm}
{\includegraphics[width=6cm,height=4cm]{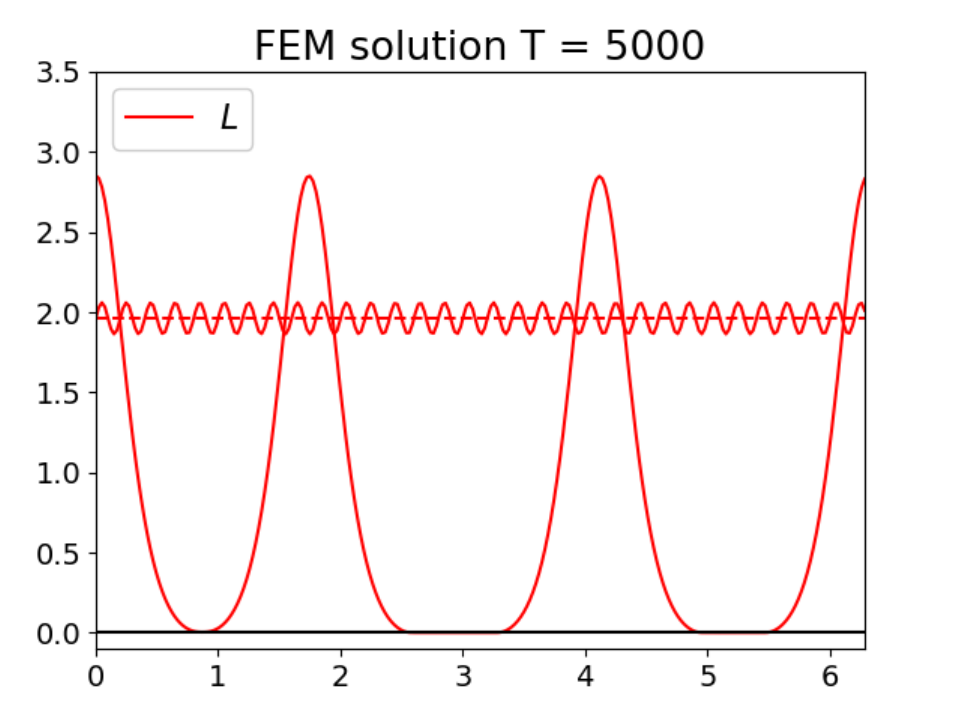}}
\caption{{\small Subcritical case. Experiment 4: Evolution of the instability formation. Only the labor is shown, the capital following a similar trend.  For large times, uninhabited regions arise separating 
densely populated areas. }} 
\label{exp3_fig}
\end{figure}

\begin{figure}[t!]
\centering
{\includegraphics[width=6cm,height=4cm]{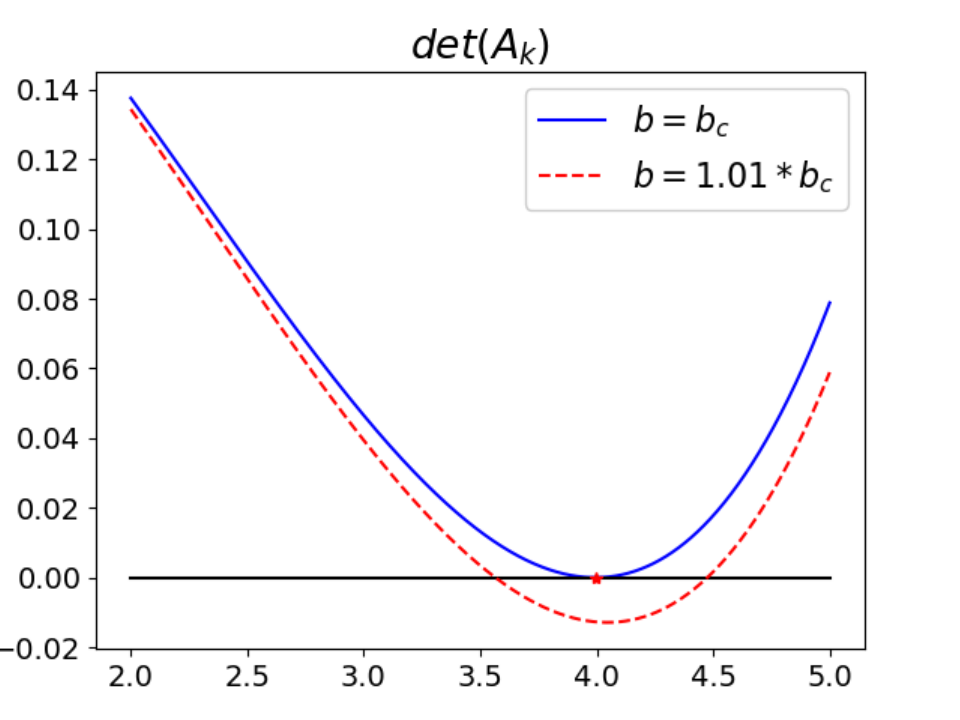}}
\hspace{0.1cm}
{\includegraphics[width=6cm,height=4cm]{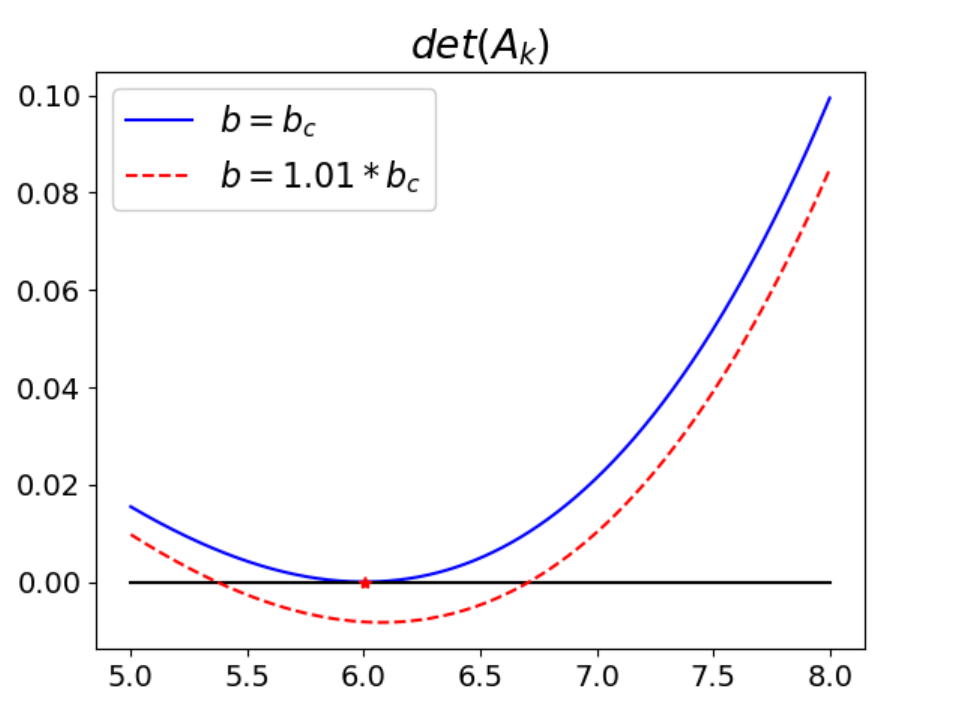}} \\
{\includegraphics[width=6cm,height=4cm]{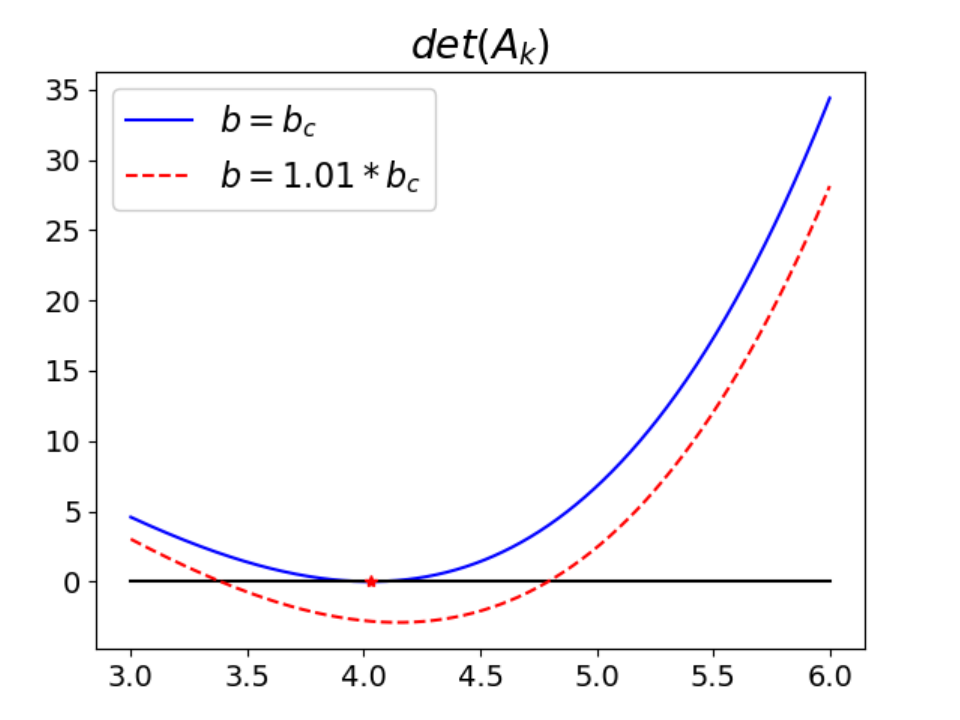}}
\hspace{0.1cm}
{\includegraphics[width=6cm,height=4cm]{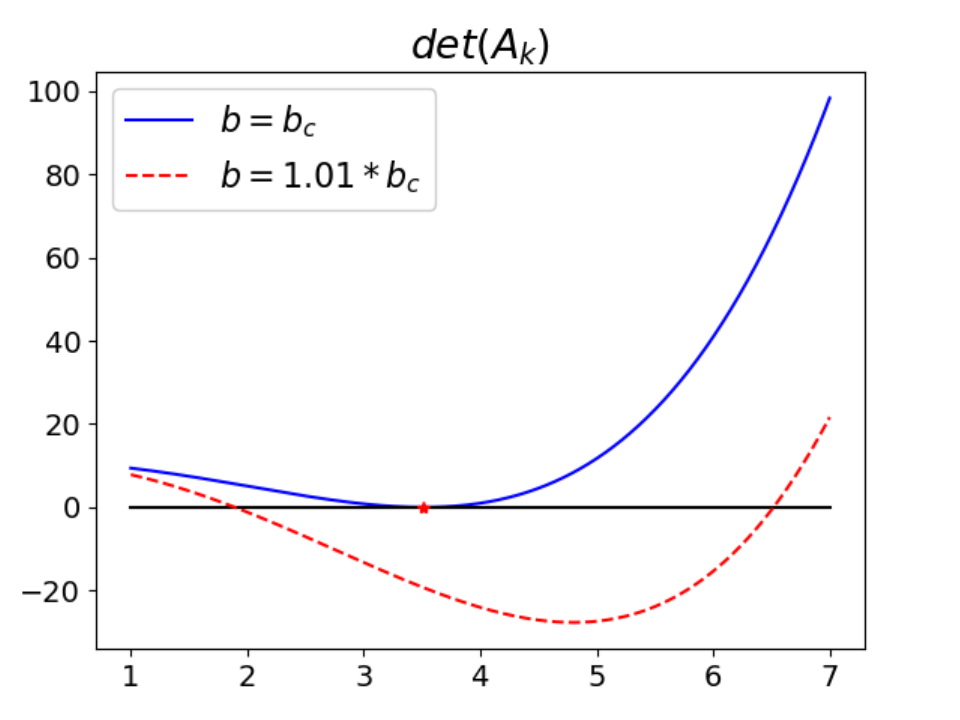}}
\caption{{\small Experiments 1 to 4: Determinant of the eigenvalue matrix $A_k$, see \fer{def:Ak}, as a function of $k$, and for two values of the bifurcation parameter: the critical value $b_c$ and the value used in the experiments, $b=1.01b_c$. 
}} 
\label{det_fig}
\end{figure}

We run four simulations with data leading to pattern formation. 
Unless otherwise stated, the following parameters are fixed for all the experiments.
The initial data is a perturbation of the stable equilibrium given by 
$L_0 = L^*(1+0.05\sin(10\pi x))$ and $K_0 = K^*(1+0.05\cos(10\pi x))$.
The time step is $\tau = 0.01$, and the spatial domain is $\O=(0,2\pi$), with 256 spatial nodes. The tolerances for the fixed point method and for the stationary state are taken as $\text{tol}_{FP} =1.e-6$  and $\text{tol}_S=1.e-8$. The total factor productivity is $A=1$, the time-scale parameter is $\gamma=1$, the saturation constant is $K_s=10K^*$, and the bifurcation parameter is $b=1.01b_c$. 
In Table~\ref{tab:data}, we show the parameters which vary in each experiment and additional information on the numerical results.
In Figures \ref{exp1_fig}-\ref{det_fig}, we show the approximated solutions obtained in the experiments and the determinant of 
the matrix $A_k$, see \fer{def:Ak}, showing the unstable wave numbers.

\begin{remark}
For developed economies, the values $\alpha\approx 0.3 $ and  $\beta\approx 0.6$ are well established. 
The values of $\epsilon$ and $\beta$ are more difficult to estimate. They can be obtained from the definition of the CES production function \fer{def:ces}, from setting the scaling variable $A$ to a known value, and  from data indicating that the ratio $K/Y \approx 3.5$ holds in a variety of situations. Of course, the choice of these parameters is not unique. In any case, we found that their role in the stability analysis is 
limited so, for the examples, we fixed them in order to obtain interesting visualizations. 
\end{remark}

In the first two experiments, the main reason because instability arises, see \fer{cond:bif} and \fer{cond:bif2},  is that the diffusion parameters promoting stability, $c_1,~c_2$, and $a_1$ are small in comparison with the ratios $\alpha_i/\det(B)$, determining the growth capacity of the system. In the Experiment 1, this happens in a symmetric way for the competition coefficients,  $\beta_1\approx\beta_2$ and in an asymmetric way for the intrinsic growth coefficients, with $\alpha_1\gg \alpha_2$. In the second experiment, the opposite relation is considered. We see that the critical bifurcation  parameter is, at least, one order higher than the other diffusion coefficients, implying that, under these situations, instability arises if the attraction felt by labor for capital is high
in comparison with random or repulsion effects.

In the third experiment we set larger stability promoting diffusion coefficients and smaller values of the growth and competition coefficients. However, it is clear that 
the important relationship for instability is that of the size of $\alpha_i/\det(B)$ against the size of the diffusion coefficients. Thus, instability emerges as well in this case. However, notice that here the critical bifurcation  parameter is of the same order as the other diffusion coefficients. This is, it is not necessary a 
extreme difference between attraction and repulsion for instability to arise:
growth plays a fundamental role.
In this experiment, we set $A = 100$ and $\gamma=5$ for a better visualization.

In the Experiment 4 we set parameters leading to large amplitude instabilities containing regions abandoned by labor and with few remaining of capital. Unlike the other experiments, in this case we experimented with a subcritical bifurcation  so that the third order Stuart-Landau equation for the amplitude does not give any useful information.  The stabilizing factor  $\beta_2 a_1$, see \fer{cond:bif}, is large, forcing a large bifurcation coefficient
for instabilities to arise. Observe that $\beta_2$ is the intra-competition coefficient for $K$, and $a_1$ is the intra-repulsion coefficient for $L$. Thus, even if these intra-population stabilizing mechanisms are intense, a large enough attraction of labor for capital is still able to cause severe instabilities. 
In this experiment, we set $ \gamma=10$.

Finally, let us notice the good agreement between the FEM and the NWL approximations in the Experiments 1 to 3, for which the bifurcation is supercritical and the third order Stuart-Landau equation provides an estimate for the steady state amplitude. In all of them, the  mean square error  between the FEM and the NWL approximations is  of the order $\text{MSE}\approx 1.e-2$.

\begin{remark}
\label{rem:saturation} 
The best place for understanding the contribution of the saturation term to instability is the sufficient condition  \fer{cond:bif}. We have that: (i) the saturation term promotes instability,  and (ii) it regulates the size needed for the bifurcation parameter to overpass the critical threshold.

In the experiments 1 and 4 this term is small in comparison to the other diffusion terms. In the experiment 2, we have $a_2\sim c_1,c_2$  and $a_2 \gg a_1$. Finally, in the experiment 3 we have $a_2 \gg a_1,c_1,c_2$ . However, instability arises in all these experiments. Therefore, it does not seem that the saturation term plays any special role for the emergence of instabilities. However, this term is crucial for the existence theory. Indeed, if this term were not uniformly bounded then the well-posedness of the problem  could not be guaranteed due to the lack of parabolicity. 
\end{remark}

\section*{Appendix: Weakly nonlinear analysis in 1D}\label{sec:wnl}

Let $\bw$ be the random perturbation around the coexistence equilibrium introduced in Section~\ref{sec:turing}. We recast the nonlinear system \fer{eq:u}-\fer{eq:v} as 
\begin{align}
\label{eq:wexp}
 \p_t \bw = \cL^b \bw + \cN^b \bw,
\end{align}
where $\cL^b = \gamma R + Q^b \p_{xx}$, and $\cN^b$ is a nonlinear operator. Here, we introduced in the notation the superscript $b$ to stress the dependence of these operators on the bifurcation parameter. We  consider the decomposition  
\begin{align*}
 \cN^b = \frac{1}{2} \big( \cQ_R(\bw,\bw) + \p_{xx} \cQ_Q(\bw,\bw) \big) + \cS_Q^b(\bw),
\end{align*}
with the bilinear symmetric forms
\begin{align*}
&  \cQ_R(\bx,\by) = \gamma \begin{pmatrix}
                            -2\beta_{1}x_1y_1 + (x_1y_2+x_2y_1) \\
                            -2\beta_{2}x_2y_2 + (x_1y_2+x_2y_1) 
                           \end{pmatrix},
\quad   \cQ_Q(\bx,\by) = \begin{pmatrix}
                            a_{1} x_1y_1\\
                            a_{2}g'(K^*)x_2y_2 
                           \end{pmatrix},
\end{align*}
and the nonlinear operator
\begin{align*}
 \cS_Q^b (\bw) = \begin{pmatrix}
                  -b \p_x(w_1 \p_x w_2 )\\
                  a_2 \sum_{j=2}^\infty \frac{g^{(j)}(K^*)}{(j+1)!}\p_{xx}w_2^{j+1}
                 \end{pmatrix}.
\end{align*}
Introducing the parameter $\eps^2 = (b-b_c)/b_c$, 
and the expansions \fer{expansion1}-\fer{expansion3},  
replacing them in the system \fer{eq:wexp}, and equating in terms of the order of $\eps$, leads to  the following systems of equations, that we set in $(0,T)\times(0,2\pi/k_c)$ and complement with non-flow boundary conditions:
\begin{align}
O(\eps): \cL^{b_c}\bw_1 =&   ~ 0,\label{prob:e1}\\
O(\eps^2): \cL^{b_c}\bw_2  =&  ~  \p_{T_1}\bw_1 - B_1\p_{xx}\bw_1 -\frac{1}{2}\big(\cQ_R(\bw_1,\bw_1) +  \p_{xx} \cQ_Q(\bw_1,\bw_1) \big)  \nonumber\\
&+ b_c \p_x (w_{1,1}\p_x w_{1,2}) \be_1 =: \bF , \label{prob:e2}\\
 O(\eps^3):  \cL^{b_c}\bw_3 = & ~ \p_{T_1}\bw_2 + \p_{T_2}\bw_1  
 - \big(B_2 \p_{xx} \bw_1+  B_1 \p_{xx} \bw_2 \big) \nonumber \\
 &  - 
\big(\cQ_R(\bw_1,\bw_2) + \p_{xx} \cQ_Q(\bw_1,\bw_2) \big) \nonumber\\
&+ b_c \p_x \big( w_{1,1}\p_x w_{2,2} + w_{2,1}\p_x w_{1,2}\big) \be_1  \nonumber\\
&+ b_1 \p_x \big( w_{1,1}\p_x w_{1,2} \big) \be_1  
+ \frac{a_2}{6}g''(K^*) \p_{xx}w_{1,2}^3  \be_2 := \bG , \label{prob:e3}
\end{align}
where $\cL^{b_c} = \gamma R + Q^{b_c} \p_{xx}$, and the elements of the matrix $B_j$
are zero, with the exception of $(B_j)_{12}= -b_j L^*$.
We compute the solutions of  \fer{prob:e1}-\fer{prob:e3}.

\no\textbf{Order $\eps$: }The solution of the linear problem \fer{prob:e1} is given by $\bw_1 = A(T_1,T_2)\brho \cos(k_cx)$, with $\brho\in \ker(\gamma R - k_c^2 Q^{b_c})$,
where $A$ is the amplitude of the pattern, unknown at the moment. Observe that $\gamma R - k_c^2 Q^{b_c} = A_{k_c}^{b_c}$, with $A_k$ defined in \fer{def:Ak}, and with $k_c$ determined to yield $\det(A_{k_c}^{b_c})=0$. Therefore, the vector $\brho$ is defined up to a multiplicative constant, that  we shall fix later, see \fer{def:rho}.

\no\textbf{Order $\eps^2$: }
We have $B_1\p_{xx} \bw_ 1 = -Ak_c^2 \cos(k_c x) B_1 \brho$.  
On noting that $\cQ_U(\bw_1,\bw_1) = A^2 \cQ_U(\brho,\brho) \cos^2(k_cx)$, for $U=R,~Q$, 
and using standard trigonometric identities, we find that 
\begin{align*}
 \frac{1}{2}\big(\cQ_R(\bw_1,\bw_1) +  \p_{xx} \cQ_Q(\bw_1,\bw_1) \big) = \frac{1}{4}A^2 \sum_{j=0,2}\cM_j(\brho,\brho) \cos(jk_cx),
\end{align*}
with $ \cM_j = \cQ_R - j^2k_c^2 \cQ_Q.$ Computing the other terms of \fer{prob:e2} yields the problem 
\begin{align*}
 \bF = &    \big( \brho  \p_{T_1} A + Ak_c^2  B_1 \brho \big) \cos(k_cx) 
  - \frac{1}{4}A^2 \sum_{j=0,2}\cM_j(\brho,\brho) \cos(jk_cx) \\
&   - A^2b_c k_c^2  \rho_1\rho_2 \be_1  \cos(2k_cx).
\end{align*}
By Fredholm's alternative, \fer{prob:e2} admits a solution if and only if 
$\langle\bF,\bpsi\rangle_{L^2} = 0$, where $\langle\cdot,\cdot\rangle_{L^2} $ denotes the scalar product in $L^2(0,2\pi/k_c)$, and $\bpsi \in \ker( (\cL^{b_c})^*)$ is of the form 
\begin{align}
\label{def:psi}
\bpsi = \boldeta \cos(k_cx), \qtext{with }\boldeta \in \ker( (\gamma R - k_c^2 Q^{b_c})^*). 
\end{align}
For similar reasons than $\brho$, $\boldeta$ is defined up to a multiplicative constant. We fix $\boldeta$ in \fer{def:rho}, where it is also shown that  $\langle \brho,\boldeta\rangle \neq 0$.

The compatibility condition implies that the terms in $\eps$ of the expansions of $b$ and $T$ are secular terms, and thus  we impose $T_1\equiv 0$ and $b_1\equiv 0$, implying $A\equiv A(T_2)$.
With these restrictions, the Fredholm's alternative is satisfied, and we look for a solution of \fer{prob:e2} of the form $ \bw_2 = A^2\sum_{j=0,2} \bw_{2j}\cos(jk_cx)$,
for which 
\begin{align*}
 \cL^{b_c} \bw_2=  A^2\sum_{j=0,2} L_j\bw_{2j}\cos(jk_cx), \qtext{with }L_j = \gamma R- j^2k_c^2Q^{b_c}.
\end{align*}
Then, $\cL^{b_c} \bw_2 =\bF$ if the vectors $\bw_{2j}$ are the solutions of the linear systems
\begin{align*}
 &L_0 \bw_{20} = -\frac{1}{4} \cM_0(\brho,\brho),\quad 
 L_2 \bw_{22} = -\frac{1}{4} \cM_j(\brho,\brho)- b_c k_c^2  \rho_1\rho_2 \be_1 .
\end{align*}


\no\textbf{Order $\eps^3$: }
Replacing $\bw_1,\bw_2$ in the terms of $\bG$, see  \fer{prob:e3},  yields
\begin{align*}
 \p_{T_2}\bw_1 = \brho \cos(k_cx) \p_{T_2}A,\qquad 
   B_2 \p_{xx} \bw_1 = -A k_c^2 B_2\brho  \cos(k_cx). 
\end{align*}
Using that $\cQ_R$ and $\cQ_Q$ are bilinear and recalling the definition of $\cM$, we get 
\begin{align*}
 \cQ_R(\bw_1,\bw_2) + \p_{xx} \cQ_Q
 (\bw_1,\bw_2)  
 =& A^3 \Big( \cos(k_cx) \big(\cM_1(\brho,\bw_{20}) +  \frac{1}{2} \cM_1(\brho,\bw_{22}) \big) \nonumber \\
 & +  \frac{1}{2} \cos(3k_cx) \cM_3(\brho,\bw_{22}) \Big). 
\end{align*}
Computing the other nonlinear terms and replacing them in the definition of $\bG$ given in \fer{prob:e3},  we obtain 
\begin{align*}
 \bG =&  \Big(\brho \p_{T_2} A   +A k_c^2 B_2\brho 
 - A^3 \big(\cM_1(\brho,\bw_{20}) +  \frac{1}{2} \cM_1(\brho,\bw_{22}) \big) \\
&  - A^3 k_c^2 b_c \rho_1 w^{(2)}_{2,2} \be_1  -\frac{1}{2} A^3 k_c^2 b_c\rho_2   \big( 2w^{(1)}_{2,0} - w^{(1)}_{2,2} \big)   \be_1\\
&   + \frac{1}{8}A^3 a_2 g''(K^*) k_c^2 \rho_2^3 \be_2 \Big) \cos(k_c x)\\
 & - A^3 \Big( \frac{1}{2} \cM_3(\brho,\bw_{22})+ 3k_c^2 b_c\rho_1 w^{(2)}_{2,2} \be_1 + \frac{3}{2}k_c^2 b_c\rho_2   w^{(1)}_{2,2} \be_1 \\
 & + \frac{3}{8}a_2 g''(K^*) k_c^2 \rho_2^3 \be_2 \Big)  \cos(3k_cx).
\end{align*}
The solvability condition for problem \fer{prob:e3} is $\langle\bG,\bpsi\rangle_{L^2} = 0$, with $\bpsi =\boldeta\cos(k_cx)$ given by \fer{def:psi}. This condition leads to the differential equation
\begin{align*}
  \langle \brho,\boldeta\rangle \p_{T_2} A + \langle\bG_1,\boldeta\rangle A +
  \langle\bG_3 ,\boldeta\rangle A^3 = 0,
\end{align*}
where $\bG_1 =  k_c^2 B_2\brho = -k_c^2L^* b_2 \rho_2 \be_1$, and 
\begin{align*}
 \bG_3 =& -\big(\cM_1(\brho,\bw_{20}) +  \frac{1}{2} \cM_1(\brho,\bw_{22}) \big) - k_c^2 b_c \rho_1w^{(2)}_{2,2} \be_1 \\
& -\frac{1}{2} k_c^2 b_c \rho_2   \big( 2w^{(1)}_{2,0} - w^{(1)}_{2,2} \big)   \be_1 
+ \frac{1}{8} a_2 g''(K^*) k_c^2 \rho_2^3 \be_2.
\end{align*}
Thus, we deduce the cubic Stuart-Landau equation for the amplitude
\begin{align}
\label{eq:sl}
 \p_{T_2} A = \sigma A - \ell A^3,
\end{align}
with
$ \sigma = -\langle \bG_1, \boldeta\rangle / \langle \brho ,\boldeta\rangle$ and 
$ \ell =   \langle \bG_3, \boldeta\rangle / \langle \brho ,\boldeta\rangle$. 
We, finally, fix the vectors $\brho\in \ker(\gamma R - k_c^2 Q^{b_c})$ , and 
$\boldeta \in \ker( (\gamma R - k_c^2 Q^{b_c})^*)$. Since $\gamma R_{21} - k_c^2 Q_{21}^{b_c} = \gamma K^* > 0$, see \fer{def:Q}, we choose the  forms 
$ \brho = (M,1)^t$, and $\boldeta =(1,M^*)^t$, where 
\begin{align}
\label{def:rho}
M= \frac{-\gamma R_{22}+ k_c^2 Q_{22}^{b_c}}{\gamma R_{21}-k_c^2Q_{21}^{b_c}}, \quad
 M^* =  \frac{-\gamma R_{11}+k_c^2 Q_{11}^{b_c}}{\gamma R_{21}-k_c^2Q_{21}^{b_c}}.
\end{align}
With this election, we have $\langle \bG_1, \boldeta\rangle = -k_c^2 b_2 L^*  < 0$, and $\langle \brho ,\boldeta\rangle >0$, 
and therefore, the growth rate coefficient $\sigma$ is always positive. As mentioned at the beginning of Section~\ref{sec:numerics}, the  dynamics of \fer{eq:sl} in the supercritical case, $\ell>0 $, has the stable equilibrium solution $A_\infty =\sqrt{\sigma/\ell}$, representing the asymptotic value of the amplitude. The corresponding solution may be approximated by \fer{wnl2}.

\end{document}